\documentclass[12pt,a4paper]{article}
\usepackage{amsfonts}
\usepackage{amssymb}
\usepackage{mathrsfs}
\usepackage{amsmath}
\usepackage{amsthm}
\usepackage{enumerate}
\usepackage{cite}
\usepackage[all]{xy}
\usepackage{extarrows}
\usepackage{indentfirst}
\allowdisplaybreaks
\newtheorem{defn}{Definition}[section]
\newtheorem{thm}[defn]{Theorem}
\newtheorem{lem}[defn]{Lemma}
\newtheorem{prop}[defn]{Proposition}
\newtheorem{cor}[defn]{Corollary}
\newtheorem{eg}[defn]{Example}
\newtheorem{re}[defn]{Remark}

\newcommand\relphantom[1]{\mathrel{\phantom{#1}}}
\newcommand{\bdefn}{\begin{defn}}
\newcommand{\edefn}{\end{defn}}
\newcommand{\bthm}{\begin{thm}}
\newcommand{\ethm}{\end{thm}}
\newcommand{\blem}{\begin{lem}}
\newcommand{\elem}{\end{lem}}
\newcommand{\bprop}{\begin{prop}}
\newcommand{\eprop}{\end{prop}}
\newcommand{\bcor}{\begin{cor}}
\newcommand{\ecor}{\end{cor}}
\newcommand{\beg}{\begin{eg}}
\newcommand{\eeg}{\end{eg}}
\newcommand{\bre}{\begin{re}}
\newcommand{\ere}{\end{re}}
\newcommand{\bpf}{\begin{proof}}
\newcommand{\epf}{\end{proof}}
\newcommand{\benu}{\begin{enumerate}}
\newcommand{\eenu}{\end{enumerate}}
\newcommand{\bc}{\begin{center}}
\newcommand{\ec}{\end{center}}
\newcommand{\bea}{\begin{eqnarray}}
\newcommand{\eea}{\end{eqnarray}}
\newcommand{\Bea}{\begin{eqnarray*}}
\newcommand{\Eea}{\end{eqnarray*}}
\newcommand{\beq}{\begin{equation}}
\newcommand{\eeq}{\end{equation}}
\newcommand{\Beq}{\begin{equation*}}
\newcommand{\Eeq}{\end{equation*}}
\newcommand{\bspl}{\begin{split}}
\newcommand{\espl}{\end{split}}

\begin{document}

\title{\textbf{On split Leibniz triple systems}
\author{ Yan Cao$^{1,2},$  Liangyun Chen$^{1}$
 \date{{\small {$^1$ School of Mathematics and Statistics, Northeast Normal
 University,\\
Changchun 130024, China}\\{\small {$^2$  Department of Basic
 Education,
 Harbin University of
Science and Technology\\ Rongcheng Campus,  Rongcheng 264300,
China}}}}}} \maketitle
\date{}

\begin{abstract}

In order to  study the structure of arbitrary split Leibniz triple
systems, we introduce the class  of split Leibniz triple systems as
the natural extension of the class of split Lie triple systems and
split Leibniz algebras. By developing techniques of connections of
roots for this kind of triple systems, we show that any of such
Leibniz triple systems $T$ with a symmetric root system is of the
form $T=U+\sum_{[j]\in \Lambda^{1}/\sim} I_{[j]}$ with $U$ a
subspace of $T_{0}$  and any $I_{[j]}$ a well described ideal of
$T$, satisfying $\{I_{[j]},T,I_{[k]}\}
=\{I_{[j]},I_{[k]},T\}=\{T,I_{[j]},I_{[k]}\}=0$ if $[j]\neq [k]$.\\

\end{abstract}
\renewcommand{\thefootnote}{\fnsymbol{footnote}}
\footnote[0]{ Corresponding author(L. Chen): chenly640@nenu.edu.cn.}
\footnote[0]{Supported by  NNSF of China (Nos. 11171055 and
11471090),  NSF of Jilin province (No. 201115006), Scientific
Research Fund of Heilongjiang Provincial Education Department
 (No. 12541184). }

\section{Introduction}

The notion of Leibniz algebras was introduced by Loday \cite{L2}, which is a ``nonantisymmetric'' generalization of Lie algebras. So far, many results of this kind of algebras have been considered in the frameworks of low dimensional
algebras, nilpotence and related problems \cite{AAO, AAO2, AO, B3, B4}.  Leibniz triple systems were introduced by Bremner and S\'{a}nchez-Ortega \cite{BS}. Leibniz triple systems were defined in a functorial manner using
 the Kolesnikov-Pozhidaev algorithm, which took the defining identities for a variety of algebras and produced the defining identities for the corresponding variety
  of dialgebras \cite{K}.  In \cite{BS}, Leibniz triple systems were obtained by applying the Kolesnikov-Pozhidaev algorithm to Lie triple systems. In \cite{BL5238999},
 Levi's theorem for Leibniz triple systems is determined.  Furthermore, Leibniz triple
  systems are related to Leibniz algebras in the same way that Lie triple systems related to Lie algebras. So it is natural to prove analogs of results from the theory
   of Lie triple systems to Leibniz triple systems.

  Recently, in \cite{BL52, BL52567, BL5234, BL528 }, the structures of arbitrary split Lie algebras, arbitrary split
  Leibniz algebras and arbitrary split  Lie triple systems have been determined by the techniques of connections of roots.
   Our work is essentially motivated by the work on  split Leibniz algebras and  split  Lie triple systems\cite{BL52567, BL5234}.

Throughout this paper, Leibniz triple systems $T$ are considered of
arbitrary dimension and over an arbitrary field $\mathbb{K}$. It is
worth to mention that, unless otherwise stated, there is not any
restriction on dim$T_{\alpha}$ or $\{k \in \mathbb{K}$:  $k \alpha
\in \Lambda^{1},$ for a fixed $\alpha \in \Lambda^{1}\}$, where
$T_{\alpha}$ denotes the root space associated to the root $\alpha$,
and $\Lambda^{1}$ the set of nonzero roots of $T$. This paper
proceeds as follows. In section 2, we establish the preliminaries on
split Leibniz triple systems theory. In section 3, we show that such
an arbitrary Leibniz triple system with a symmetric root system is
of the form  $T=U+\sum_{[j]\in \Lambda^{1}/\sim} I_{[j]}$ with $U$ a
subspace of  $T_{0}$  and any $I_{[j]}$ a well described ideal of
$T$,  satisfying $\{I_{[j]},T,I_{[k]}\}
=\{I_{[j]},I_{[k]},T\}=\{T,I_{[j]},I_{[k]}\}=0$ if $[j]\neq [k]$.

\section{Preliminaries}

\bdefn{\rm\cite{BL5234}}  A \textbf{right Leibniz algebra} $L$ is a vector space over a field $\mathbb{K}$ endowed with a bilinear product
$[\cdot,\cdot]$ satisfying the Leibniz identity
$$[[y, z], x] = [[y, x], z] + [y, [z, x]],$$
for all $x, y, z \in L$.
\edefn

\bdefn{{\rm\cite{BS}}} A \textbf{Leibniz triple system} is a vector space $T$ endowed with a trilinear
operation $\{\cdot,\cdot,\cdot\}: T\times T\times T\rightarrow T$ satisfying
\begin{gather}
\{a,\{b,c,d\},e\}\!=\!\{\{a,b,c\},d,e\}
\!-\!\{\{a,c,b\},d,e\}\!-\!\{\{a,d,b\},c,e\}\!+\!\{\{a,d,c\},b,e\},\label{VIP1}\\
\{a,b,\{c,d,e\}\}\!=\!\{\{a,b,c\},d,e\}\!-\!\{\{a,b,d\},c,e\}\!-\!\{\{a,b,e\},c,d\}\!+\!\{\{a,b,e\},d,c\},\label{VIP2}
\end{gather}
for all $a, b, c, d, e \in T$.
\edefn

\beg\label{eg}
A Lie triple system gives a Leibniz triple system with the same ternary product.
If $L$ is a Leibniz algebra with product $[\cdot,\cdot]$, then $L$ becomes a Leibniz triple system by putting $\{x,y,z\}=[[x,y],z]$. More examples refer to \cite{BS}.
\eeg

\bdefn{{\rm\cite{BS}}} Let $I$ be a subspace of a Leibniz triple system  $T$. Then $I$ is called a \textbf{subsystem} of $T$, if $\{I,I,I\}\subseteq I;$ $I$ is called an \textbf{ideal} of $T$,
if $\{I,T,T\}+\{T,I,T\}+\{T,T,I\}\subseteq I$.
\edefn

\bprop{{\rm\cite{BL5238999}}}\label{38888}
Let $T$ be a  Leibniz triple system. Then the following assertions hold.

$\rm(1)$  $J$ is generated by $\{\{a,b,c\}-\{a,c,b\}+\{b,c,a\}: a,b,c \in T\}$, then $J$ is an ideal of $T$ satisfying $\{T,T,J\}=\{T,J,T\}=0$.

$\rm(2)$  $J$ is generated by $\{\{a,b,c\}-\{a,c,b\}+\{b,c,a\}: a,b,c \in T\}$, then $T$ is a Lie triple system if and only if $J=0$.

$\rm(3)$ $\{\{c,d,e\},b,a\}-\{\{c,d,e\},a,b\}-\{\{c,b,a\},d,e\}+\{\{c,a,b\},d,e\}-\{c,\{a,b,d\},e\}-\{c,d,\{a,b,e\}\}=0$, for all $a, b, c, d, e$ $\in$ $T$.
\eprop

\bdefn{{\rm\cite{BS}}}\label{uni Leib envelop} The \textbf{standard embedding} of a Leibniz triple system $T$ is the two-graded right Leibniz algebra $L =
L^{0}\oplus L^{1}$, $L^{0}$ being the $\mathbb{K}$-$\rm span$ of $\{x \otimes y,  \ x, y \in T \}$, $L^{1}: =T$ and where the product is
given by
$$[(x \otimes y, z), (u \otimes v,w)]:= (\{x, y, u\} \otimes  v - \{x, y, v\} \otimes u
+ z \otimes w, \{x, y, w\} +\{z,u, v\}-\{z,v, u \}).$$
\edefn

Let us observe that $L^{0}$  with the product induced by the one in $L =
L^{0}\oplus L^{1}$  becomes a
 right Leibniz algebra.

\bdefn Let $T$ be a Leibniz triple system, $L = L^{0}\oplus L^{1}$ be its standard embedding, and $H^{0}$ be a  maximal abelian
subalgebra $($\rm MASA$)$
of $L^{0}$. For a linear functional $\alpha \in (H^{0})^{\ast},$
 we define the root space of $T$ $($with respect to
$H^{0}$$)$ associated to $\alpha$ as the subspace $T_{\alpha}:= \{t_{\alpha} \in T: [t_{\alpha}, h] = \alpha(h)t_{\alpha}$ for any h $\in H^{0}\}$.
The elements $\alpha \in (H^{0})^{\ast}$
 satisfying $T_{\alpha} \neq 0$ are called roots of $T$ with respect to $H^{0}$ and
we denote $\Lambda^{1}:= \{\alpha \in (H^{0})^{\ast}\setminus \{0\}: T_{\alpha}\neq 0 \}.$
\edefn

Let us observe that $T_{0} = \{t_{0} \in T: [t_{0}, h] = 0$ for any h $\in H^{0}\}$. In the following, we shall
denote by $\Lambda^{0}$ the set of all nonzero $\alpha \in (H^{0})^{\ast}$
such that $L_{\alpha}^{0}
:= \{ v_{\alpha}^{0}\in
 L^{0}: [v_{\alpha}^{0}, h]
 = \alpha(h)v_{\alpha}^{0}$
for any h $\in H^{0}\} \neq 0$.

 \blem \label{355} Let $T$ be a Leibniz triple system,  $L =L^{0}\oplus L^{1}$ be its standard embedding,
and  $H^{0}$ be a $\rm MASA$ of $L^{0}$. For $\alpha,\beta,\gamma
\in \Lambda^{1}\cup \{0\}$ and $\delta \in \Lambda^{0} \cup \{0\}$,
the following assertions hold.

$\rm(1)$ If $[T_{\alpha}, T_{\beta}]\neq 0$ then $\alpha+\beta \in  \Lambda^{0}\cup \{0\}$  and  $[T_{\alpha}, T_{\beta}]\subseteq L_{\alpha+\beta}^{0}$.

$\rm(2)$  If $[L_{\delta}^{0}, T_{\alpha}]\neq 0$ then $\delta+\alpha \in  \Lambda^{1}\cup \{0\}$  and  $[L_{\delta}^{0}, T_{\alpha}]\subseteq T_{\delta+\alpha}$.

$\rm(3)$  If $[T_{\alpha}, L_{\delta}^{0} ]\neq 0$ then $\alpha+\delta \in  \Lambda^{1}\cup \{0\}$  and  $[T_{\alpha}, L_{\delta}^{0}]\subseteq T_{\alpha+\delta}$.

$\rm(4)$ If $[L_{\delta}^{0}, L_{\gamma}^{0}]\neq 0$ then $\delta+\gamma \in  \Lambda^{0}\cup \{0\}$ and $[L_{\delta}^{0}, L_{\gamma}^{0} ]\subseteq L_{\delta+\gamma}^{0}$.

$\rm(5)$  If $\{T_{\alpha}, T_{\beta}, T_{\gamma}\}\neq 0$ then $\alpha+\beta+\gamma \in  \Lambda^{1}\cup \{0\}$  and  $\{T_{\alpha}, T_{\beta},T_{\gamma}\}\subseteq T_{\alpha+\beta+\gamma}$.
\elem

\bpf (1) For any $x \in T_{\alpha}$, $y \in T_{\beta}$ and $h \in H^{0}$, by Leibniz identity, one has $[[x,y],h]=[x,[y,h]]+[[x,h],y]=[x, \beta(h)y]+[\alpha(h)x,y]=(\alpha+\beta)(h)[x,y].$

(2)  For any $x \in L_{\delta}^{0}$, $y \in T_{\alpha}$ and $h \in H^{0}$, by Leibniz identity, one has $[[x,y],h]=[x,[y,h]]+[[x,h],y]=[x, \alpha(h)y]+[\delta(h)x,y]=(\delta+\alpha)(h)[x,y].$

(3)  For any $x \in T_{\alpha}$, $y \in L_{\delta}^{0}$,  and $h \in H^{0}$, by Leibniz identity, one has $[[x,y],h]=[x,[y,h]]+[[x,h],y]=[x, \delta(h)y]+[\alpha(h)x,y]=(\alpha+\delta)(h)[x,y].$

(4) For any $x \in L_{\delta}^{0}$, $y \in L_{\gamma}^{0}$ and $h \in H^{0}$, by Leibniz identity, one has $[[x,y],h]=[x,[y,h]]+[[x,h],y]=[x, \gamma(h)y]+[\delta(h)x,y]=(\delta+\gamma)(h)[x,y].$

(5) It is a consequence of Lemma $\ref{355}$ (1) and (2).
\epf

\bdefn Let $T$ be a Leibniz triple system, $L = L^{0}\oplus L^{1}$ be its standard embedding, and  $H^{0}$ be a $\rm MASA$
of $L^{0}$. We shall call that $T$ is a \textbf{split Leibniz triple system} $($with respect to $H^{0}$$)$ if $:$

$\rm(1)$ $T = T_{0}\oplus(\oplus_{\alpha \in \Lambda^{1}} T_{\alpha})$,

$\rm(2)$  $\{T_{0}, T_{0}, T_{0}\}=0$,

$\rm(3)$ $\{T_{\alpha}, T_{-\alpha}, T_{0}\}=0$, for $\alpha \in \Lambda^{1}$.

\noindent We say that $\Lambda^{1}$ is the root system of $T$.
\edefn

We also note that the facts $H^{0}\subset L^{0}=[T,T]$ and $T=T_{0}\oplus(\oplus_{\alpha \in \Lambda^{1}} T_{\alpha})$ imply
\beq\label{111}
H^{0}=[T_{0},T_{0}]+\sum_{\alpha \in \Lambda^{1}}[T_{\alpha},T_{-\alpha}].
\eeq

Finally, as $[T_{0}, T_{0}]\subset L^{0}_{0}=H^{0},$ we have
\beq\label{222} [T_{0},[T_{0},T_{0}]]=0. \eeq

We finally note that $\alpha \in \Lambda^{1}$ does not imply $\alpha \in \Lambda^{0}$.

\bdefn
A root system $\Lambda^{1}$ of a split Leibniz triple system $T$ is called \textbf{symmetric} if it satisfies that $\alpha \in  \Lambda^{1}$
implies $-\alpha \in \Lambda^{1}$.

A similar concept applies to the set $\Lambda^{0}$ of nonzero roots of $L^{0}$.
\edefn

In the following, $T$ denotes a split Leibniz triple system with a symmetric root system $\Lambda^{1}$, and
$T = T_{0}\oplus(\oplus_{\alpha \in \Lambda^{1}} T_{\alpha})$ the corresponding root decomposition.  We begin the study of split Leibniz triple systems by developing the concept of connections of roots.

\bdefn\label{333555}
Let $\alpha$ and $\beta$ be two nonzero roots, we shall say that $\alpha$ and $\beta$ are \textbf{connected} if there exists a family  $\{\alpha_{1},\alpha_{2},\cdots,\alpha_{2n},\alpha_{2n+1}\}\subset \Lambda^{1}\cup \{0\}$ of roots of $T$ such that

\noindent $\rm(1)$ $\{\alpha_{1},\alpha_{1}+\alpha_{2}+\alpha_{3},\alpha_{1}+\alpha_{2}+\alpha_{3}+\alpha_{4}+\alpha_{5},\cdots,\alpha_{1}+\cdots+\alpha_{2n}+\alpha_{2n+1}\} \subset \Lambda^{1}; $

\noindent $\rm(2)$ $\{\alpha_{1}+\alpha_{2}, \alpha_{1}+\alpha_{2}+\alpha_{3}+\alpha_{4},\cdots, \alpha_{1}+\cdots+\alpha_{2n}\}\subset \Lambda^{0};$

\noindent $\rm(3)$  $\alpha_{1}=\alpha$ and $\alpha_{1}+\cdots+\alpha_{2n}+\alpha_{2n+1}\in \pm\beta.$

\noindent We shall also say
that $\{\alpha_{1},\alpha_{2},\cdots,\alpha_{2n},\alpha_{2n+1}\}$ is a connection from $\alpha$ to $\beta$.
\edefn

We denote by
$$\Lambda_{\alpha}^{1}:= \{\beta \in \Lambda^{1}:  \alpha \ and \ \beta \ are \ connected\},$$
we can easily get that $\{\alpha\}$ is a connection from $\alpha$ to itself and to $-\alpha$. Therefore $\pm\alpha \in \Lambda_{\alpha}^{1}$.

\bdefn
A subset $\Omega^{1}$ of a root system $\Lambda^{1}$, associated to a split Leibniz triple system $T$, is called a\textbf{ root
subsystem} if it is symmetric, and for $\alpha, \beta, \gamma \in \Omega^{1}\cup \{0\} $ such that $\alpha + \beta \in \Lambda^{0}$ and
$\alpha +  \beta+ \gamma  \in  \Lambda^{1}$ then $\alpha +  \beta+ \gamma \in \Omega^{1}$.
\edefn

Let $\Omega^{1}$ be a root subsystem of $\Lambda^{1}$. We define
$$T_{0, \Omega^{1}} := \mathrm{span}_{\mathbb{K}}\{\{T_{\alpha}, T_{\beta}, T_{\gamma} \}:  \alpha +  \beta+ \gamma = 0;\  \alpha, \beta, \gamma \in  \Omega^{1}\cup\{0\}\} \subset T_{0}$$
and $V_{\Omega^{1}}:=\oplus_{\alpha \in \Omega^{1}}T_{\alpha}$. Taking into account the fact that $\{T_{0}, T_{0}, T_{0}\}=0$, it is straightforward to verify that
$T_{\Omega^{1}}:=T_{0,\Omega^{1}}\oplus V_{\Omega^{1}}$ is a  subsystem of $T$.  We will say that $ T_{\Omega^{1}}$ is a  subsystem associated to the
root subsystem $\Omega^{1}$.

\bprop \label{678}
If $\Lambda^{0}$ is symmetric, then the relation $\sim$ in $\Lambda^{1}$, defined by $\alpha \sim \beta$ if and only if $\beta \in \Lambda_{\alpha}^{1},$
is of equivalence.
\eprop

\bpf
This can be proved completely analogously to \cite[Proposition 3.1]{BL523}.
\epf

\bprop \label{67890}
 Let $\alpha$ be a nonzero root and suppose $\Lambda^{0}$ is symmetric. Then  $\Lambda_{\alpha}^{1}$
 is a root subsystem.
\eprop

\bpf
This can be proved completely analogously to \cite[Lemma 3.1]{BL523}.
\epf

\section{Decompositions}

 In this section, we will show that for a fixed $\alpha_{0} \in \Lambda^{1}$,
the subsystem $T_{\Lambda_{\alpha_{0}}^{1}}$ associated to the root subsystem $\Lambda_{\alpha_{0}}^{1}$ is an ideal of $T$.

 \blem \label{lemma 3.1}
The following assertions hold.

$\rm(1)$ If $\alpha, \beta \in \Lambda^{1}$ with $[T_{\alpha}, T_{\beta}] \neq 0$, then $\alpha$ is connected with $\beta.$

$\rm(2)$ If $\alpha, \beta \in \Lambda^{1}$, $\alpha \in \Lambda^{0}$  and $[L_{\alpha}^{0}, T_{\beta}] \neq 0$, then $\alpha$ is connected with $\beta.$

$\rm(3)$  If $\alpha, \beta \in \Lambda^{1}$, $\alpha \in \Lambda^{0}$  and $[T_{\beta}, L_{\alpha}^{0} ] \neq 0$, then $\alpha$ is connected with $\beta.$

$\rm(4)$ If $\alpha, \beta \in \Lambda^{1}$, $\alpha, \beta \in \Lambda^{0}$  and $[L_{\alpha}^{0}, L_{\beta}^{0}] \neq 0$, then $\alpha$  is connected with $\beta.$

$\rm(5)$ If $\alpha, \overline{\beta} \in  \Lambda^{1}$ such that $\alpha$ is not connected with $\overline{\beta}$, then $[T_{\alpha}, T _{\overline{\beta}}]= 0 $, $[L_{\alpha}^{0}, T _{\overline{\beta}}] = 0$ and $[T _{\overline{\beta}}, L_{\alpha}^{0}] = 0$ if
furthermore $\alpha  \in  \Lambda^{0} $. If $\alpha, \overline{\beta} \in  \Lambda^{1}$ such that $\alpha$ is not connected with $\overline{\beta}$,
 then $[L_{\alpha}^{0},L_{\overline{\beta}}^{0}
]= 0$  if furthermore $\alpha, \overline{\beta} \in \Lambda^{0}$.
 \elem

\bpf
(1) Suppose $[T_{\alpha}, T_{\beta}] \neq 0$, by Lemma \ref{355} (1),  one gets $\alpha+\beta \in \Lambda^{0}\cup \{0\}$. If $\alpha+\beta=0,$  then $\beta = -\alpha$ and so
$\alpha$ is connected with $\beta$. Suppose $\alpha + \beta \neq 0$. Since $\alpha+\beta \in \Lambda^{0}$, one gets $\{\alpha, \beta,-\alpha\}$ is a
connection from $\alpha$ to $\beta$.

(2) Suppose $[L_{\alpha}^{0}, T_{\beta}] \neq 0$, by Lemma \ref{355} (2),  one gets $\alpha+\beta \in \Lambda^{1}\cup \{0\}$. If $\alpha+\beta=0$, then $\beta = -\alpha$ and so
$\alpha$ is connected with $\beta$.  Suppose $\alpha + \beta \neq 0$. Since $\alpha+\beta \in \Lambda^{1}$,  we obtain $\{\alpha, 0, -\alpha-\beta\}$ is a
connection from $\alpha$ to $\beta$.

(3) Suppose $[T_{\beta}, L_{\alpha}^{0}] \neq 0$, by Lemma \ref{355} (3),  one gets $\beta+\alpha \in \Lambda^{1}\cup \{0\}$. If $\beta+\alpha=0$, then $\beta = -\alpha$ and it is clear that $\alpha$ is connected with $\beta$.  Suppose $\beta + \alpha \neq 0$. Since $\beta+\alpha \in \Lambda^{1}$,  one gets $\{\beta,  -\alpha-\beta, 0\}$ is a
connection from $\beta$ to $\alpha$. By the symmetry, we get $\alpha$ is connected with $\beta$.

(4) Suppose $[L_{\alpha}^{0}, L_{\beta}^{0}] \neq 0$, by Lemma \ref{355} (4),  one has $\alpha+\beta \in \Lambda^{0}\cup \{0\}$. If $\alpha+\beta=0$, then $\beta = -\alpha$ and so $\alpha$ is connected with $\beta$.  Suppose $\alpha + \beta \neq 0$. Since $\alpha+\beta \in \Lambda^{0}$,  one gets $\{\alpha, \beta, -\alpha\}$ is a
connection from $\alpha$ to $\beta$.

(5) It is a consequence of Lemma \ref{lemma 3.1} (1), (2),  (3) and (4).
\epf

 \blem \label{lemma 3.2}
  If $\alpha, \overline{\beta} \in  \Lambda^{1}$ are not connected, then $\{T_{\alpha}, T_{-\alpha}, T_{\overline{\beta}}\} = 0$.
\elem

\bpf
If $[T_{\alpha}, T_{-\alpha}] = 0$, it is clear. One may  suppose that $[T_{\alpha}, T_{-\alpha}] \neq 0$ and $\{T_{\alpha}, T_{-\alpha}, T_{\overline{\beta}}\}$
$\neq0$. By Leibniz identity, one gets $$\{T_{\alpha}, T_{-\alpha}, T_{\overline{\beta}}\}=[[T_{\alpha}, T_{-\alpha}], T_{\overline{\beta}}]\subset [T_{\alpha},[T_{-\alpha},T_{\overline{\beta}}]]+[[T_{\alpha},T_{\overline{\beta}}],T_{-\alpha}].$$ So either $[T_{\alpha},[T_{-\alpha},T_{\overline{\beta}}]] \neq  0$ or
$[[T_{\alpha},T_{\overline{\beta}}],T_{-\alpha}] \neq 0$, contradicting Lemma \ref {lemma 3.1} (5). Hence, $\{T_{\alpha}, T_{-\alpha}, T_{\overline{\beta}}\} = 0$.
\epf

 \blem \label{lemma 3.3}
Fix $\alpha_{0} \in \Lambda^{1}$ and suppose $\Lambda^{0}$ is
symmetric. For $\alpha \in \Lambda_{\alpha_{0}}^{1}$  and
$\beta,\gamma \in \Lambda^{1}\cup\{0\},$  the following assertions
hold.

$\rm(1)$ If $\{T_{\alpha}, T_{\beta}, T_{\gamma}\} \neq 0$  then $\beta$, $\gamma$, $\alpha+\beta+\gamma \in  \Lambda_{\alpha_{0}}^{1}
\cup \{0\}.$

$\rm(2)$ If  $\{T_{\beta}, T_{\alpha}, T_{\gamma}\} \neq 0$  then $\beta$, $\gamma$, $\beta+\alpha+\gamma \in  \Lambda_{\alpha_{0}}^{1}
\cup \{0\}.$

$\rm(3)$ If  $\{T_{\beta}, T_{\gamma}, T_{\alpha}\} \neq 0$  then $\beta$, $\gamma$, $\beta+\gamma+\alpha \in  \Lambda_{\alpha_{0}}^{1}
\cup \{0\}.$
\elem

\bpf
(1) It is easy to see that $[T_{\alpha}, T_{\beta}] \neq 0$, for $\alpha \in \Lambda_{\alpha_{0}}^{1}$  and $\beta \in \Lambda^{1}\cup\{0\}.$
By Lemma \ref{lemma 3.1} (1), one gets $\alpha \sim \beta$ in the case $\beta \neq 0$. From
here, $\beta \in \Lambda_{\alpha_{0}}^{1}\cup\{0\}$. In order to complete the proof, we will show $\gamma$, $\alpha+\beta+\gamma \in  \Lambda_{\alpha_{0}}^{1}\cup\{0\}$.
We distinguish two cases.

Case 1. Suppose $\alpha+\beta+\gamma=0$. It is clear that  $\alpha+\beta+\gamma \in \Lambda_{\alpha_{0}}^{1}\cup\{0\}$. The fact that   $\{T_{0}, T_{0}, T_{0}\}=0$ and
 $\{T_{\alpha}, T_{-\alpha}, T_{0}\}=0$ for $\alpha \in \Lambda^{1}$ gives us $\gamma \neq  0$. By Lemma \ref{355} (1),
one gets $\alpha + \beta \in \Lambda^{0}$. As $\alpha + \beta =-\gamma$,   $\{\alpha, \beta, 0\}$ would be a connection from $\alpha$ to $\gamma$ and we
conclude $\gamma \in  \Lambda_{\alpha_{0}}^{1}\cup\{0\}$.

Case 2. Suppose $\alpha+\beta+\gamma\neq 0$. We treat separately two cases.

 Suppose $\alpha + \beta \neq 0$. By  Lemma \ref{355} (1), one gets $\alpha + \beta \in \Lambda^{0}$ and so $\{\alpha, \beta, \gamma\} $ is a
connection from $\alpha$  to $\alpha+\beta+\gamma$. Hence $\alpha+\beta+\gamma \in  \Lambda_{\alpha_{0}}^{1}\cup\{0\} $. In the case $\gamma \neq 0$,
$\{\alpha, \beta, -\alpha-\beta-\gamma\}$ is a connection from $\alpha$ to  $\gamma$. So $\gamma \in \Lambda_{\alpha_{0}}^{1}$.
Hence $\gamma \in \Lambda_{\alpha_{0}}^{1}\cup \{0\}$.

 Suppose $\alpha+\beta = 0$.
Then necessarily $\gamma \in \Lambda_{\alpha_{0}}^{1}\cup\{0\}$. Indeed, if $\gamma \neq 0$ and $\alpha$ is not connected with $\gamma,$
 by Lemma \ref{lemma 3.2}, $\{T_{\alpha}, T_{\beta}, T_{\gamma}\} =\{T_{\alpha}, T_{-\alpha}, T_{\gamma}\}=0$, a contradiction.
We also have $\alpha+\beta+\gamma = \gamma \in  \Lambda_{\alpha_{0}}^{1}\cup\{0\}$.

(2) The fact that $[T_{\beta}, T_{\alpha}] \neq 0$  implies by Lemma \ref{lemma 3.1} (1) that $\alpha \sim \beta$ in the case $\beta \neq 0$. From
here, $\beta \in \Lambda_{\alpha_{0}}^{1}\cup\{0\}$.  In order to complete the proof, we will show $\gamma$, $\beta+\alpha+\gamma \in  \Lambda_{\alpha_{0}}^{1}\cup\{0\}$. We  distinguish two cases.

Case 1. Suppose $\beta+\alpha+\gamma=0$. It is clear that $\beta+\alpha+\gamma \in \Lambda_{\alpha_{0}}^{1}\cup\{0\}$. The fact that   $\{T_{0}, T_{0}, T_{0}\}=0$ and
 $\{T_{\alpha}, T_{-\alpha}, T_{0}\}=0$ for $\alpha \in \Lambda^{1}$ gives us $\gamma \neq  0$.  By Lemma \ref{355} (1),
one has  $\beta + \alpha\in \Lambda^{0}$. As $\beta + \alpha =-\gamma$,   $\{\alpha, \beta, 0\}$ would be a connection from $\alpha$ to $\gamma$ and we
conclude $\gamma \in  \Lambda_{\alpha_{0}}^{1}\cup\{0\}$.

Case 2. Suppose $\beta+\alpha+\gamma\neq 0$. We treat separately  two cases.

 Suppose $\beta + \alpha \neq 0$. By  Lemma \ref{355} (1), one gets  $\beta + \alpha \in \Lambda^{0}$ and so $\{\alpha, \beta, \gamma\} $ is a
connection from $\alpha$  to $\beta+\alpha+\gamma$. Hence $\beta+\alpha+\gamma \in  \Lambda_{\alpha_{0}}^{1}\cup\{0\} $. In the case $\gamma \neq 0$, we have
$\{\alpha, \beta, -\alpha-\beta-\gamma\}$ is a connection from $\alpha$ to  $\gamma$. So $\gamma \in \Lambda_{\alpha_{0}}^{1}$. Hence $\gamma \in \Lambda_{\alpha_{0}}^{1}\cup \{0\}$.

 Suppose $\beta+\alpha = 0$.
Then necessarily $\gamma \in \Lambda_{\alpha_{0}}^{1}\cup\{0\}$. Indeed, if $\gamma \neq 0$ and $\alpha$ is not connected with $\gamma$,
by Lemma \ref{lemma 3.2}, $\{T_{\alpha}, T_{\beta}, T_{\gamma}\} =\{T_{\alpha}, T_{-\alpha}, T_{\gamma}\}=0$, a contradiction.
We also have $\beta+\alpha+\gamma = \gamma \in  \Lambda_{\alpha_{0}}^{1}\cup\{0\}$.

(3) By Leibniz identity,  $\{T_{\beta}, T_{\gamma}, T_{\alpha}\}=[[T_{\beta}, T_{\gamma}], T_{\alpha}]\subset [T_{\beta},[T_{\gamma},T_{\alpha}]]+[[T_{\beta},T_{\alpha}],T_{\gamma}]$. From $\{T_{\beta}, T_{\gamma}, T_{\alpha}\}\neq 0$,  we  obtain  either $[T_{\beta},[T_{\gamma},T_{\alpha}]]\neq 0$ or $[[T_{\beta},T_{\alpha}],T_{\gamma}]\neq 0.$ We treat separately  two cases.

Csse 1. Suppose  $[T_{\beta},[T_{\gamma},T_{\alpha}]]\neq 0$, we will show $\beta$, $\gamma$, $\beta+\gamma+\alpha \in  \Lambda_{\alpha_{0}}^{1}
\cup \{0\}$. First to show $\gamma  \in  \Lambda_{\alpha_{0}}^{1}
\cup \{0\}$. The fact that $[T_{\gamma},T_{\alpha}]\neq 0$ implies by Lemma \ref{lemma 3.1} (1) that $\gamma \sim \alpha$ in the case $\gamma\neq 0$. From here, $
\gamma \in  \Lambda_{\alpha_{0}}^{1}\cup\{0\}$.

 Next to show $\beta \in \Lambda_{\alpha_{0}}^{1}\cup\{0\}$. Indeed, if $\beta\neq 0$ and  suppose $\alpha$ is not connected with
$\beta$, then $\beta$ is not connected with $\gamma$ in the case $\gamma\neq 0$. By Lemma \ref{lemma 3.1} (1), $[T_{\beta}, T_{\gamma}]=0$ whenever $\gamma\neq 0$, contradicting $\{T_{\beta}, T_{\gamma}, T_{\alpha}\} \neq 0$. Next to show if $\beta\neq 0$ and in the case $\gamma=0$, we  also  get  $\beta \in \Lambda_{\alpha_{0}}^{1}$. Indeed,  suppose $\alpha$ is not connected with
$\beta$, in the case $\gamma=0$, one has $\{T_{\beta}, T_{\gamma}, T_{\alpha}\} =\{T_{\beta}, T_{0}, T_{\alpha}\}=[[T_{\beta}, T_{0}], T_{\alpha}]$. From $[T_{\beta}, T_{0}]\subset L_{\beta}^{0}$ and Lemma \ref{lemma 3.1} (5), one gets  $\{T_{\beta}, T_{0}, T_{\alpha}\}=0$, a contradiction.

 Finally, to show  $\beta+\gamma+\alpha \in  \Lambda_{\alpha_{0}}^{1}\cup\{0\}$. Suppose $\beta+\gamma+\alpha=0$ and so $\beta+\gamma+\alpha \in \Lambda_{\alpha_{0}}^{1}
\cup \{0\}$. Suppose $\beta+\gamma+\alpha \neq 0$, by $[T_{\gamma},T_{\alpha}]\neq 0$,  $\gamma+\alpha \in \Lambda^{0}\cup\{0\}$. If $\gamma+\alpha \neq 0$, then $
\gamma+\alpha \in \Lambda^{0}$ and so $\{\alpha,\gamma,\beta\}$ is a connection from $\alpha$ to  $\beta+\gamma+\alpha$. Hence $\beta+\gamma+\alpha \in  \Lambda_{\alpha_{0}}^{1}$. If $\gamma+\alpha = 0$, then necessarily $\beta \in \Lambda_{\alpha_{0}}^{1}\cup\{0\}$. Indeed, if $\beta \neq 0$ and $\alpha$ is not connected with $\beta$, by Lemma \ref{lemma 3.1} (5), $\{T_{\beta}, T_{\gamma}, T_{\alpha}\} =\{T_{\beta}, T_{-\alpha}, T_{\alpha}\}=[[T_{\beta}, T_{-\alpha}], T_{\alpha}]=0$,  contradicting $\{T_{\beta}, T_{\gamma}, T_{\alpha}\} \neq 0$.
Therefore, we also have $\beta+\gamma+\alpha = \beta \in  \Lambda_{\alpha_{0}}^{1}\cup\{0\}$.

Csse 2.  If $[[T_{\beta},T_{\alpha}],T_{\gamma}]\neq 0$, we will show $\beta$, $\gamma$, $\beta+\gamma+\alpha \in  \Lambda_{\alpha_{0}}^{1}
\cup \{0\}$. First to show $
\beta \in  \Lambda_{\alpha_{0}}^{1}\cup\{0\}$. The fact that $[T_{\beta},T_{\alpha}]\neq 0$ implies by Lemma \ref{lemma 3.1} (1) that $\beta \sim \alpha$ in the case $\beta \neq 0$. From here, $
\beta \in  \Lambda_{\alpha_{0}}^{1}\cup\{0\}$.

 Next to show $\gamma \in \Lambda_{\alpha_{0}}^{1}\cup\{0\}$.
Indeed, if $\gamma\neq 0$ and $\alpha$ is not connected with
$\gamma$, then $\beta$ is not connected with $\gamma$ in the case $\beta\neq 0$. By Lemma \ref{lemma 3.1} (1), $[T_{\beta}, T_{\gamma}]=0$ whenever $\beta\neq 0$, contradicting $\{T_{\beta}, T_{\gamma}, T_{\alpha}\} \neq 0$. Similarly, it is easy to show if $\gamma\neq 0$ and in the case $\beta=0$, we can obtain  $\gamma \in \Lambda_{\alpha_{0}}^{1}$.

 Finally, to show  $\beta+\gamma+\alpha \in  \Lambda_{\alpha_{0}}^{1}\cup\{0\}$. Suppose $\beta+\gamma+\alpha=0$ and so $\beta+\gamma+\alpha \in \Lambda_{\alpha_{0}}^{1}
\cup \{0\}$. Suppose $\beta+\gamma+\alpha \neq 0$, by
$[T_{\beta},T_{\alpha}]\neq 0$,  one has  $\beta+\alpha \in
\Lambda^{0}\cup\{0\}$. If $\beta+\alpha \neq 0$, then $ \beta+\alpha
\in \Lambda^{0}$ and so $\{\alpha,\beta,\gamma\}$ is a connection
from $\alpha$ to  $\beta+\gamma+\alpha$. Hence $\beta+\gamma+\alpha
\in  \Lambda_{\alpha_{0}}^{1}$. If $\beta+\alpha = 0$, then
necessarily $\gamma \in \Lambda_{\alpha_{0}}^{1}\cup\{0\}$. Indeed,
if $\gamma \neq 0$ and $\alpha$ is not connected with $\gamma$, by
Lemma \ref{lemma 3.1} (5), $\{T_{\beta}, T_{\gamma}, T_{\alpha}\}
=\{T_{-\alpha}, T_{\gamma}, T_{\alpha}\}=[[T_{-\alpha}, T_{\gamma}],
T_{\alpha}]=0$,  contradicting $\{T_{\beta}, T_{\gamma},
T_{\alpha}\} \neq 0$. Therefore, we also have $\beta+\gamma+\alpha =
\gamma \in  \Lambda_{\alpha_{0}}^{1}\cup\{0\}$.\epf

\blem \label{lemma 3.4}
Fix $\alpha_{0} \in \Lambda^{1}$ and suppose $\Lambda^{0}$ is symmetric. For $\alpha,\beta,\gamma \in \Lambda_{\alpha_{0}}^{1}\cup\{0\}$ with $\alpha+\beta+\gamma=0$ and $\delta,\epsilon \in \Lambda^{1}\cup\{0\},$
 the following assertions hold.

$\rm(1)$ If $\{\{T_{\alpha}, T_{\beta}, T_{\gamma}\},T_{\delta}, T_{\epsilon}\}\neq 0$  then $\delta$, $\epsilon$, $\delta+\epsilon \in  \Lambda_{\alpha_{0}}^{1}
\cup \{0\}.$

$\rm(2)$ If $\{T_{\delta}, \{T_{\alpha}, T_{\beta}, T_{\gamma}\}, T_{\epsilon}\}\neq 0$  then $\delta$, $\epsilon$, $\delta+\epsilon \in  \Lambda_{\alpha_{0}}^{1}
\cup \{0\}.$

$\rm(3)$ If $\{T_{\delta}, T_{\epsilon}, \{T_{\alpha}, T_{\beta}, T_{\gamma}\} \}\neq 0$  then $\delta$, $\epsilon$, $\delta+\epsilon \in  \Lambda_{\alpha_{0}}^{1}
\cup \{0\}.$

\elem

\bpf
(1) From the fact that $\alpha + \beta + \gamma = 0$ , $\{T_{0}, T_{0}, T_{0}\}=0$ and $\{T_{\alpha}, T_{-\alpha}, T_{0}\}=0$ whenever $\alpha \in \Lambda^{1}$, one may suppose that at least two distinct
elements in $\{\alpha, \beta, \gamma\}$ are nonzero and one may  consider the case  $\{T_{\alpha}, T_{\beta}, T_{\gamma} \} \neq 0$, $\alpha + \beta \neq 0$
and $ \gamma \neq 0$. By (\ref{VIP2}), one gets
\begin{align*}
0\neq  \{\{T_{\alpha}, T_{\beta}, T_{\gamma}\},T_{\delta},T_{\epsilon}\} &\subset \{T_{\alpha}, T_{\beta},\{ T_{\gamma},T_{\delta},T_{\epsilon}\}\}+\{\{T_{\alpha}, T_{\beta}, T_{\delta}\},T_{\gamma},T_{\epsilon}\}\\
&+\{\{T_{\alpha}, T_{\beta}, T_{\epsilon}\},T_{\gamma},T_{\delta}\}+\{\{T_{\alpha}, T_{\beta}, T_{\epsilon}\},T_{\delta},T_{\gamma}\},
\end{align*}
\noindent any of the above four summands is nonzero. In order to prove $\delta$, $\epsilon$, $\delta+\epsilon \in  \Lambda_{\alpha_{0}}^{1}
\cup \{0\},$ we will consider four cases.

Case 1. Suppose $\{T_{\alpha}, T_{\beta},\{ T_{\gamma},T_{\delta},T_{\epsilon}\}\}\neq 0.$
 As $\gamma\neq 0$ and $\{ T_{\gamma},T_{\delta},T_{\epsilon}\}\neq 0$, Lemma \ref{lemma 3.3} (1) shows that $\delta, \epsilon, \gamma+\delta+\epsilon $ are
connected with $\gamma$ in the case of being nonzero roots and so $\delta, \epsilon, \gamma+\delta+\epsilon  \in \Lambda_{\alpha_{0}}^{1}\cup\{0\}$.
If $\gamma+\delta+\epsilon=0$, then $\delta+\epsilon=-\gamma \in \Lambda_{\alpha_{0}}^{1}.$ If $\gamma+\delta+\epsilon \neq 0$, taking into account $0\neq
\{T_{\alpha}, T_{\beta},\{ T_{\gamma},T_{\delta},T_{\epsilon}\}\}\subset \{T_{\alpha}, T_{\beta}, T_{\gamma+\delta+\epsilon}\}$, Lemma \ref{lemma 3.3} (3) gives
 us that $\alpha+\beta+\gamma+\delta+\epsilon=\delta+\epsilon \in \Lambda_{\alpha_{0}}^{1}$.

 Case 2. Suppose $\{\{T_{\alpha}, T_{\beta}, T_{\delta}\},T_{\gamma},T_{\epsilon}\} \neq 0$. It is clear that  $\{T_{\alpha}, T_{\beta}, T_{\delta}\}\neq 0$. As $\alpha+\beta \neq 0$, one gets either $\alpha \in \Lambda_{\alpha_{0}}^{1}$ or $\beta \in \Lambda_{\alpha_{0}}^{1}$. By Lemma \ref{lemma 3.3} (1) and (2), one gets $\delta \in \Lambda_{\alpha_{0}}^{1}
\cup \{0\}$. It is obvious that  $0\neq \{\{T_{\alpha}, T_{\beta}, T_{\delta}\},T_{\gamma},T_{\epsilon}\} \subset \{T_{\alpha+\beta+\delta}, T_{\gamma},  T_{\epsilon}\}$. As $\gamma\neq 0$, $\gamma \in \Lambda_{\alpha_{0}}^{1}$, by  Lemma \ref{lemma 3.3} (2), one gets $\epsilon \in \Lambda_{\alpha_{0}}^{1}\cup\{0\}$ and $\alpha+\beta+\gamma+\delta+\epsilon=\delta+\epsilon \in  \Lambda_{\alpha_{0}}^{1}\cup\{0\}$.

Case 3. Suppose $\{\{T_{\alpha}, T_{\beta}, T_{\epsilon}\},T_{\gamma},T_{\delta}\}\neq 0$. It is easy to see that  $\{T_{\alpha}, T_{\beta}, T_{\epsilon}\}\neq 0$. As $\alpha+\beta \neq 0$, we get either $\alpha \in \Lambda_{\alpha_{0}}^{1}$ or $\beta \in \Lambda_{\alpha_{0}}^{1}$. By Lemma \ref{lemma 3.3} (1) and (2), one gets $\epsilon \in \Lambda_{\alpha_{0}}^{1}
\cup \{0\}$. Note that  $0\neq \{\{T_{\alpha}, T_{\beta}, T_{\epsilon}\},T_{\gamma},T_{\delta}\} \subset  \{T_{\alpha+\beta+\epsilon}, T_{\gamma},  T_{\delta}\}$. As $\gamma\neq 0$, $\gamma \in \Lambda_{\alpha_{0}}^{1}$, by  Lemma \ref{lemma 3.3} (2), one gets $\delta \in \Lambda_{\alpha_{0}}^{1}\cup\{0\}$ and $\alpha+\beta+\gamma+\delta+\epsilon=\delta+\epsilon \in  \Lambda_{\alpha_{0}}^{1}\cup\{0\}$.

Case 4. Suppose $\{\{T_{\alpha}, T_{\beta}, T_{\epsilon}\},T_{\delta},T_{\gamma}\}\neq 0$. It is clear that $\{T_{\alpha}, T_{\beta}, T_{\epsilon}\}\neq 0$. As $\alpha+\beta \neq 0$, one gets either $\alpha \in \Lambda_{\alpha_{0}}^{1}$ or $\beta \in \Lambda_{\alpha_{0}}^{1}$. By Lemma \ref{lemma 3.3} (1) and (2), one gets $\epsilon \in \Lambda_{\alpha_{0}}^{1}
\cup \{0\}$. It is clear that $0\neq \{\{T_{\alpha}, T_{\beta}, T_{\epsilon}\},T_{\delta},T_{\gamma}\} \subset  \{T_{\alpha+\beta+\epsilon}, T_{\delta},  T_{\gamma}\}$. As $\gamma\neq 0$, $\gamma \in \Lambda_{\alpha_{0}}^{1}$, by  Lemma \ref{lemma 3.3} (3), one gets $\delta \in \Lambda_{\alpha_{0}}^{1}\cup\{0\}$ and $\alpha+\beta+\epsilon+\delta+\gamma=\delta+\epsilon \in  \Lambda_{\alpha_{0}}^{1}\cup\{0\}$.

(2) By Proposition \ref{38888} (3), we obtain that
 \begin{align*}
  0&\neq \{T_{\delta}, \{T_{\alpha}, T_{\beta},T_{\gamma}\},T_{\epsilon}\}\subset \{\{T_{\delta}, T_{\gamma}, T_{\epsilon}\},T_{\beta},T_{\alpha}\}+
\{\{T_{\delta}, T_{\gamma}, T_{\epsilon}\},T_{\alpha},T_{\beta}\}\\
  & +\{\{T_{\delta}, T_{\beta}, T_{\alpha}\},T_{\gamma},T_{\epsilon}\}+\{\{T_{\delta}, T_{\alpha}, T_{\beta}\},T_{\gamma},T_{\epsilon}\}+\{T_{\delta}, T_{\gamma}, \{T_{\alpha},T_{\beta},T_{\epsilon}\}\},
 \end{align*}
\noindent any of the above five summands is nonzero.

 Suppose $\{\{T_{\delta}, T_{\gamma}, T_{\epsilon}\},T_{\beta},T_{\alpha}\}\neq 0$, it is obvious $\{T_{\delta}, T_{\gamma}, T_{\epsilon}\}\neq 0$. As $\gamma\neq 0$, $\gamma \in \Lambda_{\alpha_{0}}^{1}$, by  Lemma \ref{lemma 3.3} (2), one gets $\delta, \delta+\gamma+\epsilon, \epsilon
 \in  \Lambda_{\alpha_{0}}^{1}\cup\{0\}$. Note that $0\neq \{\{T_{\delta}, T_{\gamma}, T_{\epsilon}\},T_{\beta},T_{\alpha}\} \subset  \{T_{\delta+\gamma+\epsilon}, T_{\beta},  T_{\alpha}\}$. As $\alpha+\beta \neq 0$, we get either $\alpha \in \Lambda_{\alpha_{0}}^{1}$ or $\beta \in \Lambda_{\alpha_{0}}^{1}$. By Lemma \ref{lemma 3.3} (2) and (3), one gets $\delta+\gamma+\epsilon+\beta+\alpha=\delta+\epsilon \in \Lambda_{\alpha_{0}}^{1}
\cup \{0\}$.

  If $\{\{T_{\delta}, T_{\gamma}, T_{\epsilon}\},T_{\alpha},T_{\beta}\}\neq 0$, $\{\{T_{\delta}, T_{\beta}, T_{\alpha}\},T_{\gamma},T_{\epsilon}\}\neq 0$,
$\{\{T_{\delta}, T_{\alpha}, T_{\beta}\},T_{\gamma},T_{\epsilon}\}\neq 0$ or $\{T_{\delta}, T_{\gamma}, \{T_{\alpha},T_{\beta},T_{\epsilon}\}\}\neq 0$,  a similar argument gives us $\delta$, $\epsilon$, $\delta+\epsilon \in  \Lambda_{\alpha_{0}}^{1}
\cup \{0\}.$

(3) By Proposition \ref{38888} (3), we obtain that
\begin{align*}
0&\neq\{T_{\delta}, T_{\epsilon}, \{T_{\alpha},T_{\beta},T_{\gamma}\}\}\subset \{\{T_{\delta}, T_{\epsilon}, T_{\gamma}\},T_{\beta},T_{\alpha}\}+
\{\{T_{\delta}, T_{\epsilon}, T_{\gamma}\},T_{\alpha},T_{\beta}\}\\
 &+\{\{T_{\delta}, T_{\beta}, T_{\alpha}\},T_{\epsilon},T_{\gamma}\}+\{\{T_{\delta}, T_{\alpha}, T_{\beta}\},T_{\epsilon},T_{\gamma}\}+\{T_{\delta}, \{T_{\alpha}, T_{\beta},T_{\epsilon}\}, T_{\gamma}\},
 \end{align*}
 any of the above five summands is nonzero.

 Suppose $\{\{T_{\delta}, T_{\epsilon}, T_{\gamma}\},T_{\beta},T_{\alpha}\}\neq 0$.  One easily gets $\{T_{\delta}, T_{\epsilon}, T_{\gamma}\}\neq 0$. As $\gamma\neq 0$, $\gamma \in \Lambda_{\alpha_{0}}^{1}$, by  Lemma \ref{lemma 3.3} (3), one has $\delta, \delta+\epsilon+\gamma, \epsilon
 \in  \Lambda_{\alpha_{0}}^{1}\cup\{0\}$. Note that $0\neq \{\{T_{\delta}, T_{\epsilon}, T_{\gamma}\},T_{\beta},T_{\alpha}\} \subset  \{T_{\delta+\epsilon+\gamma}, T_{\beta},  T_{\alpha}\}$. As $\alpha+\beta \neq 0$, one  gets either $\alpha \in \Lambda_{\alpha_{0}}^{1}$ or $\beta \in \Lambda_{\alpha_{0}}^{1}$. By Lemma \ref{lemma 3.3} (2) and (3), one has $\delta+\epsilon+\gamma+\beta+\alpha=\delta+\epsilon \in \Lambda_{\alpha_{0}}^{1}
\cup \{0\}$.

  If $\{\{T_{\delta}, T_{\epsilon}, T_{\gamma}\},T_{\alpha},T_{\beta}\}\neq 0$, $\{\{T_{\delta}, T_{\beta}, T_{\alpha}\},T_{\epsilon},T_{\gamma}\}\neq 0,$ $\{\{T_{\delta}, T_{\alpha}, T_{\beta}\},T_{\epsilon},T_{\gamma}\}\neq 0$ or $\{T_{\delta}, \{T_{\alpha}, T_{\beta},T_{\epsilon}\}, T_{\gamma}\}\neq 0$,  a similar argument gives us    $\delta$, $\epsilon$, $\delta+\epsilon \in  \Lambda_{\alpha_{0}}^{1}
\cup \{0\}.$
\epf

\blem \label{lemma 3.5}
Fix $\alpha_{0} \in \Lambda^{1}$and suppose $\Lambda^{0}$ is symmetric. If $\alpha_{1}, \alpha_{2}, \alpha_{3} \in \Lambda_{\alpha_{0}}^{1}\cup\{0\}$ with
$\alpha_{1} + \alpha_{2} + \alpha_{3} = 0$ and $\overline{\epsilon} \in \Lambda^{1} \setminus \Lambda_{\alpha_{0}}^{1}$, then the following assertions hold.

$\rm(1)$ $[\{T_{\alpha_{1}}, T_{\alpha_{2}}, T_{\alpha_{3}}\}, T_{\overline{\epsilon}} ] = 0.$

$\rm(2)$ In case $\overline{\epsilon} \in \Lambda^{0}$, then $[\{T_{\alpha_{1}}, T_{\alpha_{2}}, T_{\alpha_{3}}\}, L_{\overline{\epsilon}}^{0}]= 0. $

$\rm(3)$ $[[\{T_{\alpha_{1}},T_{\alpha_{2}},T_{\alpha_{3}}\},T_{0}], T_{\overline{\epsilon}}]=0.$

\elem

\bpf
(1) From the fact $\alpha_{1} + \alpha_{2} + \alpha_{3} = 0$, $\{T_{0}, T_{0}, T_{0}\}=0$ and $\{T_{\alpha}, T_{-\alpha}, T_{0}\}=0$ for $\alpha \in \Lambda^{1}$,   one gets  if $\alpha_{3}= 0$ then it is clear that $[\{T_{\alpha_{1}}, T_{\alpha_{2}}, T_{\alpha_{3}}\}, T_{\overline{\epsilon}} ] = 0.$
Let us consider the case $\alpha_{3} \neq 0$. By Leibniz identity, we have
\beq\label{444}
[\{T_{\alpha_{1}}, T_{\alpha_{2}}, T_{\alpha_{3}}\}, T_{\overline{\epsilon}} ]
 \!=\![[[T_{\alpha_{1}}, T_{\alpha_{2}}], T_{\alpha_{3}}], T_{\overline{\epsilon}} ]
  \!\subset\!
[[T_{\alpha_{1}}, T_{\alpha_{2}}], [T_{\alpha_{3}},T_{\overline{\epsilon}} ]]\!+\! [[[T_{\alpha_{1}}, T_{\alpha_{2}}], T_{\overline{\epsilon}}], T_{\alpha_{3}}]\!.
\eeq
Let us consider the first summand in (\ref{444}). As $\alpha_{3} \neq 0$,  one has $\alpha_{3} \in \Lambda_{\alpha_{0}}^{1}$.  For $\overline{\epsilon} \in \Lambda^{1} \setminus \Lambda_{\alpha_{0}}^{1}$ and
Lemma \ref{lemma 3.1} (5), one easily gets $[T_{\alpha_{3}},T_{\overline{\epsilon}} ]=0$. Therefore $[[T_{\alpha_{1}}, T_{\alpha_{2}}], [T_{\alpha_{3}},T_{\overline{\epsilon}} ]]=0.$

 Let us now consider the second summand    in (\ref{444}), it is sufficient to verify that  $$[[[T_{\alpha_{1}}, T_{\alpha_{2}}], T_{\overline{\epsilon}}], T_{\alpha_{3}}]=0.$$ To do so,  we first assert that $[[T_{\alpha_{1}}, T_{\alpha_{2}}], T_{\overline{\epsilon}}]=0$.  Indeed, by  Leibniz identity, we have \begin{equation}\label{8888}
 [[T_{\alpha_{1}}, T_{\alpha_{2}}], T_{\overline{\epsilon}}]\subset  [T_{\alpha_{1}}, [T_{\alpha_{2}}, T_{\overline{\epsilon}}]]+ [[T_{\alpha_{1}}, T_{\overline{\epsilon}}] ,T_{\alpha_{2}}],
 \end{equation}
 where $\alpha_{1}, \alpha_{2} \in \Lambda_{\alpha_{0}}^{1}\cup\{0\}$, $\overline{\epsilon} \in \Lambda^{1} \setminus \Lambda_{\alpha_{0}}^{1}$. In the following, we distinguish
 three cases.

 Case 1. $\alpha_{1}\neq 0$ and $\alpha_{2}\neq 0.$ As $\alpha_{1} \in \Lambda_{\alpha_{0}}^{1}$ and  $\overline{\epsilon} \in \Lambda^{1} \setminus \Lambda_{\alpha_{0}}^{1}$, by Lemma \ref{lemma 3.1} (1), one gets $[T_{\alpha_{1}},  T_{\overline{\epsilon}}]=0$.
As $\alpha_{2} \in \Lambda_{\alpha_{0}}^{1}$ and  $\overline{\epsilon} \in \Lambda^{1} \setminus \Lambda_{\alpha_{0}}^{1}$, by Lemma \ref{lemma 3.1} (1), one gets $[T_{\alpha_{2}},  T_{\overline{\epsilon}}]=0$. Therefore by (\ref{8888}),  one can show that $[[T_{\alpha_{1}}, T_{\alpha_{2}}], T_{\overline{\epsilon}}]=0$.

  Case 2. $\alpha_{1}\neq 0$ and $\alpha_{2}=0.$ As $\alpha_{1} \in \Lambda_{\alpha_{0}}^{1}$ and  $\overline{\epsilon} \in \Lambda^{1} \setminus \Lambda_{\alpha_{0}}^{1}$, by Lemma \ref{lemma 3.1} (1), one gets $[T_{\alpha_{1}},  T_{\overline{\epsilon}}]=0$.  That is $[[T_{\alpha_{1}}, T_{\overline{\epsilon}}] ,T_{\alpha_{2}}]=0$. As $\alpha_{2}=0$, $[T_{\alpha_{2}},  T_{\overline{\epsilon}}]=[T_{0},  T_{\overline{\epsilon}}]\subset L_{\overline{\epsilon}}^{0}$. By Lemma \ref{lemma 3.1} (5), one gets $[T_{\alpha_{1}}, [T_{\alpha_{2}}, T_{\overline{\epsilon}}]]=0$.  Therefore by (\ref{8888}),  one can show that  $[[T_{\alpha_{1}}, T_{\alpha_{2}}], T_{\overline{\epsilon}}]=0$.

 Case 3. $\alpha_{1}= 0$ and $\alpha_{2}\neq 0.$ As $\alpha_{2} \in \Lambda_{\alpha_{0}}^{1}$ and  $\overline{\epsilon} \in \Lambda^{1} \setminus \Lambda_{\alpha_{0}}^{1}$, by Lemma \ref{lemma 3.1} (1), one gets $[T_{\alpha_{2}},  T_{\overline{\epsilon}}]=0$. That is $[T_{\alpha_{1}}, [T_{\alpha_{2}}, T_{\overline{\epsilon}}]]=0.$
  As $\alpha_{1}=0$, $[T_{\alpha_{1}},  T_{\overline{\epsilon}}]=[T_{0},  T_{\overline{\epsilon}}]\subset L_{\overline{\epsilon}}^{0}$. By Lemma \ref{lemma 3.1} (5), we get $[ [T_{\alpha_{1}}, T_{\overline{\epsilon}}],T_{\alpha_{2}}]=0$.  Therefore by (\ref{8888}),  one can show that $[[T_{\alpha_{1}}, T_{\alpha_{2}}], T_{\overline{\epsilon}}]=0$.

So $[[[T_{\alpha_{1}}, T_{\alpha_{2}}], T_{\overline{\epsilon}}], T_{\alpha_{3}}]$=0 is a consequence of $[[T_{\alpha_{1}}, T_{\alpha_{2}}], T_{\overline{\epsilon}}]=0$. By (\ref{444}), one gets $[\{T_{\alpha_{1}}, T_{\alpha_{2}}, T_{\alpha_{3}}\}, T_{\overline{\epsilon}} ]=0$. The proof is complete.

(2)  From the fact $\alpha_{1} + \alpha_{2} + \alpha_{3} = 0$, $\{T_{0}, T_{0}, T_{0}\}=0$ and $\{T_{\alpha}, T_{-\alpha}, T_{0}\}=0$ for $\alpha \in \Lambda^{1}$,   one gets  if $\alpha_{3}= 0$ then it is clear that $[\{T_{\alpha_{1}}, T_{\alpha_{2}}, T_{\alpha_{3}}\}, L_{\overline{\epsilon}}^{0}  ]  = 0.$
Let us consider the case $\alpha_{3} \neq 0$. Note that
\beq\label{555666779999}
[\{T_{\alpha_{1}}, T_{\alpha_{2}}, T_{\alpha_{3}}\}, L_{\overline{\epsilon}}^{0} ]
 \subset
[[T_{\alpha_{1}}, T_{\alpha_{2}}], [T_{\alpha_{3}},L_{\overline{\epsilon}}^{0} ]]+[[[T_{\alpha_{1}}, T_{\alpha_{2}}], L_{\overline{\epsilon}}^{0}], T_{\alpha_{3}}].
\eeq
Let us consider the first summand in (\ref{555666779999}). As $\alpha_{3} \neq 0$,  one  gets $[[T_{\alpha_{1}}, T_{\alpha_{2}}], [T_{\alpha_{3}},L_{\overline{\epsilon}}^{0} ]]=0$ by
Lemma \ref{lemma 3.1} (5). Let us now consider the second summand in (\ref{555666779999}). As either $\alpha_{1}\neq 0$ or $\alpha_{2}\neq 0$,
Leibniz identity, the fact $[T_{0} , L_{\overline{\epsilon}}^{0}] \subset T_{\overline{\epsilon}} $ and Lemma \ref{lemma 3.1} (5), we obtain that $[[[T_{\alpha_{1}}, T_{\alpha_{2}}], L_{\overline{\epsilon}}^{0}], T_{\alpha_{3}}] = 0$. So, the second summand in (\ref{555666779999}) is also zero and then $[\{T_{\alpha_{1}}, T_{\alpha_{2}}, T_{\alpha_{3}}\}, L_{\overline{\epsilon}}^{0} ]$ = 0.

(3) It is  a consequence of Lemma \ref{lemma 3.5} (1), (2) and
$$[[\{T_{\alpha_{1}},T_{\alpha_{2}},T_{\alpha_{3}}\},T_{0}], T_{\overline{\epsilon}}]\subset [\{T_{\alpha_{1}},T_{\alpha_{2}},T_{\alpha_{3}}\},[T_{0},T_{\overline{\epsilon}}]]+[[\{T_{\alpha_{1}},T_{\alpha_{2}},T_{\alpha_{3}}\},T_{\overline{\epsilon}}], T_{0}].$$
\epf

\bdefn\label{ 3.16789}
A Leibniz triple system $T$ is said to be \textbf{simple} if its product is nonzreo and its only ideals are $\{0\}$, $J$ and $T$.
\edefn

It should be noted that the above definition agrees with the definition of a simple Lie triple system, since $J=\{0\}$ in this case.

\bthm \label{theorem 3.1} Suppose $\Lambda^{0}$ is symmetric, the following assertions hold.

$\rm(1)$ For any $\alpha_{0} \in \Lambda^{1}$, the  subsystem $$T_{\Lambda_{\alpha_{0}}^{1}}=T_{0,\Lambda_{\alpha_{0}}^{1}}\oplus V_{\Lambda_{\alpha_{0}}^{1}}$$
of $T$ associated to the root subsystem $\Lambda_{\alpha_{0}}^{1}$ is an ideal of $T$.

$\rm(2)$ If $T$ is simple, then there exists a connection from $\alpha$ to $\beta$ for any $\alpha$, $\beta \in \Lambda^{1}.$
\ethm
\bpf (1) Recall that $$T_{0, \Lambda_{\alpha_{0}}^{1}} := \mathrm{span}_{\mathbb{K}}\{\{T_{\alpha}, T_{\beta}, T_{\gamma} \}:  \alpha +  \beta+ \gamma = 0;\  \alpha, \beta, \gamma \in \Lambda_{\alpha_{0}}^{1} \cup\{0\}\} \subset T_{0}$$
and $V_{\Lambda_{\alpha_{0}}^{1}}:=\oplus_{\gamma \in \Lambda_{\alpha_{0}}^{1}}T_{\gamma}$. In order to complete the proof, it is sufficient to show that $$\{T_{\Lambda_{\alpha_{0}}^{1}},T,T\}+\{T, T_{\Lambda_{\alpha_{0}}^{1}},T\}+\{T,T, T_{\Lambda_{\alpha_{0}}^{1}}\}\subset T_{\Lambda_{\alpha_{0}}^{1}}.$$ We first check that $\{T_{\Lambda_{\alpha_{0}}^{1}},T,T\}\subset T_{\Lambda_{\alpha_{0}}^{1}}.$
It is easy to see that $$\{T_{\Lambda_{\alpha_{0}}^{1}},T,T\}=\{T_{0,\Lambda_{\alpha_{0}}^{1}}\oplus V_{\Lambda_{\alpha_{0}}^{1}},T,T\}=\{T_{0,\Lambda_{\alpha_{0}}^{1}},T,T\}+\{ V_{\Lambda_{\alpha_{0}}^{1}},T,T\}.$$
Next, we will show that $\{T_{0,\Lambda_{\alpha_{0}}^{1}},T,T\} \subset T_{\Lambda_{\alpha_{0}}^{1}}.$ Note that
\begin{align*}
\{T_{0,\Lambda_{\alpha_{0}}^{1}},T,T\}=&\{T_{0,\Lambda_{\alpha_{0}}^{1}}, T_{0}\oplus(\oplus_{\alpha \in \Lambda^{1}} T_{\alpha}), T_{0}\oplus(\oplus_{\alpha \in \Lambda^{1}} T_{\alpha})\}\\
=&\{T_{0,\Lambda_{\alpha_{0}}^{1}},T_{0},T_{0}\}+\{T_{0,\Lambda_{\alpha_{0}}^{1}},T_{0},\oplus_{\alpha \in \Lambda^{1}}T_{\alpha}\}\\
+&\{T_{0,\Lambda_{\alpha_{0}}^{1}},\oplus_{\alpha \in \Lambda^{1}}T_{\alpha},T_{0}\}
+\{T_{0,\Lambda_{\alpha_{0}}^{1}},\oplus_{\alpha \in \Lambda^{1}}T_{\alpha},\oplus_{\beta \in \Lambda^{1}}T_{\beta}\}.
\end{align*}
Here,  it is clear that $\{T_{0,\Lambda_{\alpha_{0}}^{1}},T_{0},T_{0}\}\subset \{T_{0},T_{0},T_{0}\}=0$.
Taking into account $\{T_{0,\Lambda_{\alpha_{0}}^{1}},T_{0},T_{\alpha}\}$, for $\alpha \in \Lambda^{1}$,  Lemma \ref{lemma 3.4} (1) and the fact that either $\alpha \in \Lambda_{\alpha_{0}}^{1}$ or
$\alpha \not \in \Lambda_{\alpha_{0}}^{1}$, give us that $\{T_{0,\Lambda_{\alpha_{0}}^{1}},T_{0},T_{\alpha}\}\subset V_{\Lambda_{\alpha_{0}}^{1}}$ or
$\{T_{0,\Lambda_{\alpha_{0}}^{1}},T_{0},T_{\alpha}\}=0$. Similarly,  one gets that $\{T_{0,\Lambda_{\alpha_{0}}^{1}},T_{\alpha},T_{0}\}\subset V_{\Lambda_{\alpha_{0}}^{1}}$ or
$\{T_{0,\Lambda_{\alpha_{0}}^{1}},T_{\alpha},T_{0}\}=0$. Next, we will consider  $\{T_{0,\Lambda_{\alpha_{0}}^{1}},T_{\alpha},T_{\beta}\}$, where $\alpha,\beta \in \Lambda^{1}$. We treat five cases.

 Case 1. If  $\alpha \in \Lambda_{\alpha_{0}}^{1}$, $\beta \in \Lambda_{\alpha_{0}}^{1}$ and $\alpha+\beta=0$. One has $$\{T_{0,\Lambda_{\alpha_{0}}^{1}},T_{\alpha},T_{\beta}\}\subset T_{0,\Lambda_{\alpha_{0}}^{1}}.$$

  Case 2. If $\alpha \in \Lambda_{\alpha_{0}}^{1}$, $\beta \in \Lambda_{\alpha_{0}}^{1}$ and $\alpha+\beta\neq 0$. By $\Lambda_{\alpha_{0}}^{1}$ is a root subsystem, one gets $$\{T_{0,\Lambda_{\alpha_{0}}^{1}},T_{\alpha},T_{\beta}\}\subset V_{\Lambda_{\alpha_{0}}^{1}}.$$

  Case 3.  If $\alpha \in \Lambda_{\alpha_{0}}^{1}$ and  $\beta \not \in \Lambda_{\alpha_{0}}^{1}$. By Lemma \ref{lemma 3.4} (1), one has $$\{T_{0,\Lambda_{\alpha_{0}}^{1}},T_{\alpha},T_{\beta}\}=0.$$

  Case 4.  If $\beta \in \Lambda_{\alpha_{0}}^{1}$ and  $\alpha \not \in \Lambda_{\alpha_{0}}^{1}$. By Lemma \ref{lemma 3.4} (1), one has $$\{T_{0,\Lambda_{\alpha_{0}}^{1}},T_{\alpha},T_{\beta}\}=0.$$

  Case 5.  If $\beta \not \in \Lambda_{\alpha_{0}}^{1}$ and  $\alpha \not \in \Lambda_{\alpha_{0}}^{1}$. By Lemma \ref{lemma 3.4} (1), one has $$\{T_{0,\Lambda_{\alpha_{0}}^{1}},T_{\alpha},T_{\beta}\}=0. $$

\noindent  Therefore, $\{T_{0,\Lambda_{\alpha_{0}}^{1}},T,T\}\subset T_{\Lambda_{\alpha_{0}}^{1}}.$

Next, we will show that  $\{ V_{\Lambda_{\alpha_{0}}^{1}},T,T\}\subset T_{\Lambda_{\alpha_{0}}^{1}}$. It is obvious that
 \begin{align*}
 \{ V_{\Lambda_{\alpha_{0}}^{1}},T,T\}&=\{\oplus _{\gamma \in \Lambda_{\alpha_{0}}^{1} }T_{\gamma},  T_{0}\oplus(\oplus_{\alpha \in \Lambda^{1}} T_{\alpha}), T_{0}\oplus(\oplus_{\alpha \in \Lambda^{1}} T_{\alpha})\}\\
&=\{\oplus _{\gamma \in \Lambda_{\alpha_{0}}^{1} }T_{\gamma},T_{0}, T_{0}\}+\{\oplus _{\gamma \in \Lambda_{\alpha_{0}}^{1} }T_{\gamma},T_{0},\oplus_{\alpha \in \Lambda^{1}} T_{\alpha}\}\\
&+\{\oplus _{\gamma \in \Lambda_{\alpha_{0}}^{1} }T_{\gamma},\oplus_{\alpha \in \Lambda^{1}} T_{\alpha},T_{0}\}+\{\oplus _{\gamma \in \Lambda_{\alpha_{0}}^{1} }T_{\gamma},\oplus_{\alpha \in \Lambda^{1}} T_{\alpha},\oplus_{\beta \in \Lambda^{1}} T_{\beta}\}.
\end{align*}

\noindent Here, it is clear that $\{T_{\gamma},T_{0}, T_{0}\}\subset  V_{\Lambda_{\alpha_{0}}^{1}} $, for $\gamma \in \Lambda_{\alpha_{0}}^{1}$.
 Next, we will consider $\{T_{\gamma},T_{0}, T_{\alpha}\}$, for $\gamma \in \Lambda_{\alpha_{0}}^{1}$, $\alpha \in \Lambda^{1}$.
 We treat three cases.

 Case 1. If $\gamma \in \Lambda_{\alpha_{0}}^{1}$, $\alpha \not \in \Lambda_{\alpha_{0}}^{1}$. By Lemma \ref{lemma 3.3} (1),  one has
 $$\{T_{\gamma},T_{0}, T_{\alpha}\}=0.$$

 Case 2. If $\gamma \in \Lambda_{\alpha_{0}}^{1}$, $\alpha  \in \Lambda_{\alpha_{0}}^{1}$ and $\gamma+\alpha \neq 0$. By  $\Lambda_{\alpha_{0}}^{1}$ is a root subsystem,  one has
$$\{T_{\gamma},T_{0}, T_{\alpha}\}\subset V_{\Lambda_{\alpha_{0}}^{1}}.$$

 Case 3. If $\gamma \in \Lambda_{\alpha_{0}}^{1}$, $\alpha  \in \Lambda_{\alpha_{0}}^{1}$ and $\gamma+\alpha = 0$.   It is clear that
$$\{T_{\gamma},T_{0}, T_{\alpha}\}\subset T_{0,\Lambda_{\alpha_{0}}^{1}}.$$

\noindent Hence,  $\{T_{\gamma},T_{0}, T_{\alpha}\}\subset T_{\Lambda_{\alpha_{0}}^{1}}$, for $\gamma \in \Lambda_{\alpha_{0}}^{1}$, $\alpha \in \Lambda^{1}$. Similarly, it is easy to get $\{T_{\gamma},T_{\alpha},T_{0}\}\subset T_{\Lambda_{\alpha_{0}}^{1}}, $ for $\gamma \in \Lambda_{\alpha_{0}}^{1}$, $\alpha \in \Lambda^{1}$.
At last,  we will consider $ \{\oplus _{\gamma \in \Lambda_{\alpha_{0}}^{1} }T_{\gamma},\oplus_{\alpha \in \Lambda^{1}} T_{\alpha},\oplus_{\beta \in \Lambda^{1}} T_{\beta}\},$ for $\gamma \in \Lambda_{\alpha_{0}}^{1},$ $\alpha \in \Lambda^{1}$ and $\beta \in \Lambda^{1}$. We treat five cases.

Case 1. If $\gamma \in \Lambda_{\alpha_{0}}^{1}$, $\alpha \in \Lambda_{\alpha_{0}}^{1}$, $\beta \in \Lambda_{\alpha_{0}}^{1}$ and $\gamma+\alpha+\beta=0$, one gets  $$ \{T_{\gamma}, T_{\alpha}, T_{\beta}\}\subset T_{0,\Lambda_{\alpha_{0}}^{1}}.$$

Case 2. If $\gamma \in \Lambda_{\alpha_{0}}^{1}$, $\alpha \in \Lambda_{\alpha_{0}}^{1}$, $\beta \in \Lambda_{\alpha_{0}}^{1}$ and $\gamma+\alpha+\beta\neq 0$, one gets  $$ \{\oplus _{\gamma \in \Lambda_{\alpha_{0}}^{1} }T_{\gamma},\oplus_{\alpha \in \Lambda^{1}} T_{\alpha},\oplus_{\beta \in \Lambda^{1}} T_{\beta}\}\subset V_{\Lambda_{\alpha_{0}}^{1}}.$$

Case 3. If $\gamma \in \Lambda_{\alpha_{0}}^{1}$, $\alpha \in \Lambda_{\alpha_{0}}^{1}$ and $\beta \not \in \Lambda_{\alpha_{0}}^{1}$. By Lemma \ref{lemma 3.3} (1) and (2), one gets $$ \{T_{\gamma}, T_{\alpha}, T_{\beta}\}=0.$$

 Case 4. If $\gamma \in \Lambda_{\alpha_{0}}^{1}$, $\alpha \not \in \Lambda_{\alpha_{0}}^{1}$ and $\beta  \in \Lambda_{\alpha_{0}}^{1}$.  By Lemma \ref{lemma 3.3} (1) and (3), one gets $$ \{T_{\gamma}, T_{\alpha}, T_{\beta}\}=0.$$

  Case 5. If $\gamma \in \Lambda_{\alpha_{0}}^{1}$, $\alpha \not \in \Lambda_{\alpha_{0}}^{1}$ and $\beta  \not \in \Lambda_{\alpha_{0}}^{1}$.  By Lemma \ref{lemma 3.3} (1), one gets $$ \{T_{\gamma}, T_{\alpha}, T_{\beta}\}=0.$$

\noindent So, $\{ V_{\Lambda_{\alpha_{0}}^{1}},T,T\}\subset T_{\Lambda_{\alpha_{0}}^{1}}.$ Therefore $\{T_{\Lambda_{\alpha_{0}}^{1}},T,T\}\subset T_{\Lambda_{\alpha_{0}}^{1}}$ is   a consequence of $\{T_{0,\Lambda_{\alpha_{0}}^{1}},T,T\}\subset T_{\Lambda_{\alpha_{0}}^{1}}$ and $\{ V_{\Lambda_{\alpha_{0}}^{1}},T,T\}\subset T_{\Lambda_{\alpha_{0}}^{1}}$.

 A similar argument gives us $\{T, T_{\Lambda_{\alpha_{0}}^{1}},T\}\subset T_{\Lambda_{\alpha_{0}}^{1}}$ and $\{T,T, T_{\Lambda_{\alpha_{0}}^{1}}\} \subset T_{\Lambda_{\alpha_{0}}^{1}}$. Consequently, this proves that $T_{\Lambda_{\alpha_{0}}^{1}}$ is an ideal of $T$.

 (2) The simplicity of $T$ implies $T_{\Lambda_{\alpha_{0}}^{1}} \in \{J,T\}$ for any $a \in \Lambda^{1}$. If $\alpha \in \Lambda^{1}$ is such that $T_{\Lambda_{\alpha_{0}}^{1}}=T$. Then $\Lambda_{\alpha_{0}}^{1}=\Lambda^{1}$. Hence, $T$ has all its nonzero roots connected. Otherwise, if $T_{\Lambda_{\alpha_{0}}^{1}}=J$ for any $\alpha \in \Lambda^{1}$ then $\Lambda_{\alpha_{0}}^{1}=\Lambda_{\beta_{0}}^{1}$ for any $\alpha_{0}$, $\beta_{0} \in  \Lambda^{1}$ and so $\Lambda_{\alpha_{0}}^{1}=\Lambda^{1}$. We also conclude that $T$ has all its nonzero roots connected.
 \epf

\bthm \label{theorem 3.2}
 Suppose $\Lambda^{0}$ is symmetric. Then for a vector space complement $U$ of
$\rm span$$_{\mathbb{K}}\{\{T_{\alpha}, T_{\beta}, T_{\gamma}\}: \alpha + \beta + \gamma = 0$, where $\alpha, \beta, \gamma \in \Lambda^{1}\cup \{0\}\}$ in $T_{0}$, we have
$$T=U+\sum_{[\alpha]\in \Lambda^{1}/\sim} I_{[\alpha]},$$
where any $I_{[\alpha]}$ is one of the ideals $T_{\Lambda_{\alpha_{0}}^{1}}$ of $T$ described in Theorem \ref{theorem 3.1}. Moreover \\ $\{I_{[\alpha]},T,I_{[\beta]}\}
=\{I_{[\alpha]},I_{[\beta]},T\}=\{T,I_{[\alpha]},I_{[\beta]}\}=0$ if $[\alpha]\neq [\beta]$.

\ethm

\bpf
Let us denote
$\xi_{0}:$=span$_{\mathbb{K}}\{\{T_{\alpha}, T_{\beta}, T_{\gamma}\}: \alpha + \beta + \gamma = 0$, where $\alpha, \beta, \gamma \in \Lambda^{1}\cup \{0\}\}$ in $T_{0}$. By Proposition \ref{678}, we can consider the quotient set $\Lambda^{1}/\sim:= \{[\alpha]: \alpha \in \Lambda^{1}\}$. By denoting $I_{[\alpha]}:=T_{\Lambda_{\alpha}^{1}}$, $T_{0,[\alpha]}:=T_{0,\Lambda_{\alpha}^{1}}$ and $V_{[\alpha]}:=V_{\Lambda_{\alpha}^{1}}$, one gets $I_{[\alpha]}:=T_{0,[\alpha]}\oplus V_{[\alpha]}$.
From
$$T=T_{0}\oplus (\oplus_{\alpha \in \Lambda^{1}}T_{\alpha} )=(  U+ \xi_{0})\oplus ( \oplus_{\alpha \in \Lambda^{1}}T_{\alpha} ),$$
 it follows

$\oplus_{\alpha \in \Lambda^{1}}T_{\alpha}=\oplus_{[\alpha] \in \Lambda^{1}/\sim }V_{[\alpha]}$, \quad  $\xi_{0}=\sum_{[\alpha] \in \Lambda^{1} /\sim }T_{0,[\alpha]}$,

\noindent which implies
$$T= U+ \xi_{0}\oplus ( \oplus_{\alpha \in \Lambda^{1}}T_{\alpha} )=U+\sum\limits_{[\alpha] \in \Lambda^{1}/\sim}I_{[\alpha]},$$
\noindent where each $I_{[\alpha]}$ is an ideal of $T$ by Theorem \ref{theorem 3.1}.

Next, it is sufficient to show that $\{I_{[\alpha]},T,I_{[\beta]}\}=0$ if $[\alpha]\neq [\beta]$. Note that,
 \begin{align*}
\{I_{[\alpha]},T,I_{[\beta]}\}=&\{T_{0,[\alpha]}\oplus V_{[\alpha]}, T_{0}\oplus ( \oplus_{\gamma \in \Lambda^{1}}T_{\gamma} ), T_{0,[\beta]}\oplus V_{[\beta]}\}\\
=& \{T_{0,[\alpha]},T_{0},T_{0,[\beta]}\}+\{T_{0,[\alpha]},T_{0},V_{[\beta]}\}+\{T_{0,[\alpha]},\oplus_{\gamma \in \Lambda^{1}}T_{\gamma},T_{0,[\beta]}\}\\
+&\{T_{0,[\alpha]},\oplus_{\gamma \in \Lambda^{1}}T_{\gamma},V_{[\beta]}\}+\{V_{[\alpha]},T_{0},T_{0,[\beta]}\}+\{V_{[\alpha]},T_{0},V_{[\beta]}\}\\
+&\{V_{[\alpha]},\oplus_{\gamma \in \Lambda^{1}}T_{\gamma},T_{0,[\beta]}\}+\{V_{[\alpha]},\oplus_{\gamma \in \Lambda^{1}}T_{\gamma},V_{[\beta]}\}.
\end{align*}

\noindent Here,  it is clear that $\{T_{0,[\alpha]},T_{0},T_{0,[\beta]}\} \subset \{T_{0},T_{0},T_{0}\}=0$. If $[\alpha]\neq [\beta],$ by Lemma \ref{lemma 3.3} and   \ref{lemma 3.4},  it is easy to see $\{T_{0,[\alpha]},T_{0},V_{[\beta]}\}=0,$ $\{T_{0,[\alpha]},\oplus_{\gamma \in \Lambda^{1}}T_{\gamma},V_{[\beta]}\}=0,$
 $\{V_{[\alpha]},T_{0},T_{0,[\beta]}\}=0,$ $\{V_{[\alpha]},T_{0},V_{[\beta]}\}=0,$ $\{V_{[\alpha]},\oplus_{\gamma \in \Lambda^{1}}T_{\gamma},T_{0,[\beta]}\}=0,$ $\{V_{[\alpha]},\oplus_{\gamma \in \Lambda^{1}}T_{\gamma},V_{[\beta]}\}=0.$

  Next, we will show  $\{T_{0,[\alpha]},\oplus_{\gamma \in \Lambda^{1}}T_{\gamma},T_{0,[\beta]}\}=0$. Indeed, for $\{T_{\alpha_{1}}, T_{\alpha_{2}},  T_{\alpha_{3}}\} \in T_{0,[\alpha]}$ with $\alpha_{1},\alpha_{2},\alpha_{3} \in \Lambda_{\alpha}^{1}\cup \{0\}$, $\alpha_{1}+\alpha_{2}+\alpha_{3}=0$, and  for  $\{T_{\beta_{1}}, T_{\beta_{2}},  T_{\beta_{3}}\} \in T_{0,[\beta]}$ with $\beta_{1},\beta_{2},\beta_{3} \in \Lambda_{\beta}^{1}\cup \{0\}$, $\beta_{1}+\beta_{2}+\beta_{3}=0$, by Proposition \ref{38888} (3), one gets
 \begin{align*}
&\{\{T_{\alpha_{1}}, T_{\alpha_{2}}, T_{\alpha_{3}}\}, \oplus_{\gamma \in \Lambda^{1}}T_{\gamma}, \{T_{\beta_{1}}, T_{\beta_{2}}, T_{\beta_{3}}\}\}\\
\subset &  \{\{\{T_{\alpha_{1}}, T_{\alpha_{2}}, T_{\alpha_{3}}\}, \oplus_{\gamma \in \Lambda^{1}}T_{\gamma}, T_{\beta_{3}}\}, T_{\beta_{2}}, T_{\beta_{1}}\}+ \{\{\{T_{\alpha_{1}}, T_{\alpha_{2}}, T_{\alpha_{3}}\}, \oplus_{\gamma \in \Lambda^{1}}T_{\gamma}, T_{\beta_{3}}\}, T_{\beta_{1}}, T_{\beta_{2}}\}\\
+&  \{\{\{T_{\alpha_{1}}, T_{\alpha_{2}}, T_{\alpha_{3}}\},  T_{\beta_{2}}, T_{\beta_{1}}\}, \oplus_{\gamma \in \Lambda^{1}}T_{\gamma}, T_{\beta_{3}}\}+\{\{\{T_{\alpha_{1}}, T_{\alpha_{2}}, T_{\alpha_{3}}\},  T_{\beta_{1}}, T_{\beta_{2}}\}, \oplus_{\gamma \in \Lambda^{1}}T_{\gamma}, T_{\beta_{3}}\}\\
 +& \{\{T_{\alpha_{1}}, T_{\alpha_{2}}, T_{\alpha_{3}}\}, \{ T_{\beta_{1}}, T_{\beta_{2}}, \oplus_{\gamma \in \Lambda^{1}}T_{\gamma}\}, T_{\beta_{3}}\}.
\end{align*}
By Lemma \ref{lemma 3.4}, it is easy to see that $$\{\{\{T_{\alpha_{1}}, T_{\alpha_{2}}, T_{\alpha_{3}}\}, \oplus_{\gamma \in \Lambda^{1}}T_{\gamma}, T_{\beta_{3}}\}, T_{\beta_{2}}, T_{\beta_{1}}\}=0,$$
 $$\{\{\{T_{\alpha_{1}}, T_{\alpha_{2}}, T_{\alpha_{3}}\}, \oplus_{\gamma \in \Lambda^{1}}T_{\gamma}, T_{\beta_{3}}\}, T_{\beta_{1}}, T_{\beta_{2}}\}=0,$$  $$\{\{\{T_{\alpha_{1}}, T_{\alpha_{2}}, T_{\alpha_{3}}\},  T_{\beta_{2}}, T_{\beta_{1}}\}, \oplus_{\gamma \in \Lambda^{1}}T_{\gamma}, T_{\beta_{3}}\}=0,$$
$$\{\{\{T_{\alpha_{1}}, T_{\alpha_{2}}, T_{\alpha_{3}}\},  T_{\beta_{1}}, T_{\beta_{2}}\}, \oplus_{\gamma \in \Lambda^{1}}T_{\gamma}, T_{\beta_{3}}\}=0,$$
 $$\{\{T_{\alpha_{1}}, T_{\alpha_{2}}, T_{\alpha_{3}}\}, \{ T_{\beta_{1}}, T_{\beta_{2}}, \oplus_{\gamma \in \Lambda^{1}}T_{\gamma}\}, T_{\beta_{3}}\}=0.$$
for $\alpha_{1},\alpha_{2},\alpha_{3} \in \Lambda_{\alpha}^{1}\cup \{0\}$, $\alpha_{1}+\alpha_{2}+\alpha_{3}=0$, $\beta_{1},\beta_{2},\beta_{3} \in \Lambda_{\beta}^{1}\cup \{0\}$, $\beta_{1}+\beta_{2}+\beta_{3}=0$, $[\alpha]\neq [\beta]$. So $\{I_{[\alpha]},T,I_{[\beta]}\}=0$  if $[\alpha]\neq [\beta]$ is  proved.

 A similar argument gives us $\{I_{[\alpha]},I_{[\beta]},T\}=\{T,I_{[\alpha]},I_{[\beta]}\}=0$ if $[\alpha]\neq [\beta].$
\epf

\bdefn
The \textbf{annihilator} of a Leibniz triple system $T$ is the set $\mathrm{Ann}(T)=\{x \in T: \{x, T, T\}+ \{T, x,  T\}+\{T,T,x\}=0\}$.
\edefn

\bcor
Suppose $\Lambda^{0}$ is symmetric. If $\mathrm{Ann}(T)=0$, and $\{T,T,T\}=T$, then $T$ is the direct sum of the ideals given in Theorem \ref{theorem 3.2},
$$T=\oplus_{[\alpha] \in \Lambda^{1}/\sim}I_{[\alpha]}.$$
\ecor

\bpf
From $\{T,T,T\}=T$ and  Theorem \ref{theorem 3.2}, we have $$\{ U+\sum_{[\alpha]\in \Lambda^{1}/\sim} I_{[\alpha]},U+\sum_{[\alpha]\in \Lambda^{1}/\sim} I_{[\alpha]},U+\sum_{[\alpha]\in \Lambda^{1}/\sim} I_{[\alpha]} \}=U+\sum_{[\alpha]\in \Lambda^{1}/\sim} I_{[\alpha]}.$$ Taking into account $U\subset T_{0}$,  Lemma \ref{lemma 3.3} and the fact that $\{I_{[\alpha]},T,I_{[\beta]}\}=\{I_{[\alpha]},I_{[\beta]},T\}=\{T,I_{[\alpha]},I_{[\beta]}\}=0$ if $[\alpha]\neq [\beta]$ (see Theorem \ref{theorem 3.2}) give us that $U=0$. That is, $$T=\sum_{[\alpha] \in \Lambda^{1}/\sim}I_{[\alpha]}.$$

To finish, it is sufficient to show the direct character of the sum.  For $x\in I_{[\alpha]}\cap \sum_{[\beta] \in \Lambda^{1}/\sim \atop \beta \not \sim \alpha} I_{[\beta]}$,  using again the equation $\{I_{[\alpha]},T,I_{[\beta]}\}=0$ for $[\alpha]\neq [\beta],$  we obtain $$\{x, T,  I_{[\alpha]}\}=\{x, T, \sum_{[\beta] \in \Lambda^{1}/\sim \atop \beta \not \sim \alpha} I_{[\beta]} \}=0.$$
 So $\{x, T, T\}=\{x,T, I_{[\alpha]}+\sum_{[\beta] \in \Lambda^{1}/\sim \atop \beta \not \sim \alpha} I_{[\beta]} \}=\{x,T, I_{[\alpha]}\}+\{x,T, \sum_{[\beta] \in \Lambda^{1}/\sim \atop \beta \not \sim \alpha} I_{[\beta]}\}=0+0=0.$
We argue similarly. Using the equations $\{T,I_{[\alpha]},I_{[\beta]}\}=0$ and  $\{I_{[\alpha]},I_{[\beta]},T\}=0$ for $[\alpha]\neq [\beta],$ one gets $\{T, x,  T\}=0$ and  $\{T, T, x \}=0$.
That is, $x \in \mathrm{Ann}(T)=0$. Thus $x=0$, as desired.
\epf

\end{document}